\documentclass[12pt,reqno]{amsart}

\usepackage{amssymb,latexsym, amsmath, amscd, array}

\numberwithin{equation}{section}
\numberwithin{figure}{section}
\numberwithin{table}{section}

\DeclareMathOperator{\syscat}{{\rm cat_{sys}}}
\DeclareMathOperator{\cat}{{\mbox{\rm cat$_{\rm LS}$}}}

\DeclareMathOperator{\length}{{\rm length}}
\DeclareMathOperator{\vol}{{\rm vol}}

\DeclareMathOperator{\area}{{\rm area}}

\DeclareMathOperator{\corank}{{\rm corank}}

\DeclareMathOperator{\sys}{{\rm sys} \pi_1}

\DeclareMathOperator{\CL}{{\it CL}} 

\DeclareMathOperator{\CA}{{\it CA}}

\DeclareMathOperator{\Reeb}{{\rm Reeb}}

\DeclareMathOperator{\diam}{{\rm diam}}

\DeclareMathOperator{\SR}{{\rm SR}}

\def\ov{\overline} \def\wh{\widehat}

\DeclareMathOperator{\FIG}{{\rm FIG}}

\def\rp{\mathbb R\mathbb P}

  \def\ov{\overline}

\newcommand{\N}{\mathbb N}

\newcommand{\R}{\mathbb R}

\newcommand{\Z}{\mathbb Z}
\newcommand{\gmetric}{{\mathcal G}}

\newtheorem{theorem}{Theorem}[section]

\newtheorem{proposition}[theorem]{Proposition}
\newtheorem{lemma}[theorem]{Lemma}
\newtheorem{corollary}[theorem]{Corollary}

\newtheorem{prop}[theorem]{Proposition}

\newtheorem{cory}[theorem]{Corollary}

\theoremstyle{definition}
\newtheorem{definition}[theorem]{Definition}
\newtheorem{example}[theorem]{Example}
\newtheorem{remark}[theorem]{Remark}

\newtheorem{question}[theorem]{Question}

\newcommand\theoref{Theorem~\ref}
\newcommand\lemref{Lemma~\ref}
\newcommand\propref{Proposition~\ref}
\newcommand\corref{Corollary~\ref}
\newcommand\defref{Definition~\ref}
\def\coryref{\corref}

\def\ie {{\it i.e.\ }}  
\def\cf {\hbox{\it cf.\ }}

\def\dist{{\rm dist}}

\def\eps{\varepsilon}
\def\gf{\varphi}

\long\def\forget#1\forgotten{} %

\date{\today}
\begin{document}

\author[M.~Katz]{Mikhail G. Katz$^{*}$} \address{Department of
Mathematics, Bar Ilan University, Ramat Gan 52900 Israel}
\email{katzmik@math.biu.ac.il} \thanks{$^{*}$Supported by the Israel
Science Foundation (grants no.~84/03 and 1294/06)}

\author[Y.~Rudyak]{Yuli B. Rudyak$^{**}$} \address{Department of
Mathematics, University of Florida, PO Box 118105, Gainesville, FL
32611-8105 USA} \email{rudyak@math.ufl.edu} \thanks{$^{**}$Supported
by the NSF, grant 0406311}

\author[S.~Sabourau]{St\'ephane Sabourau} \address{Laboratoire de
Math\'ematiques et Physique Th\'eorique, Universit\'e de Tours, Parc
de Grandmont, 37400 Tours, France} \address{
Department of Mathematics, University of Pennsylvania, 209 South 33rd Street, 
Philadelphia, PA 19104-6395, USA} \email{sabourau@lmpt.univ-tours.fr}

\title [Systoles of 2-complexes via Grushko decomposition] {Systoles of
2-complexes, Reeb graph, and Grushko decomposition}



\begin{abstract}
Let~$X$ be a finite 2-complex{} with unfree fundamental group.  We prove
lower bounds for the area of a metric on~$X$, in terms of the square
of the least length of a noncontractible loop in~$X$.  We thus
establish a uniform systolic inequality for all unfree~$2$-complexes.
Our inequality improves the constant in M.~Gromov's inequality in this
dimension.  The argument relies on the Reeb graph and the coarea
formula, combined with an induction on the number of freely
indecomposable factors in Grushko's decomposition of the fundamental
group.  More specifically, we construct a kind of a Reeb space
``minimal model'' for~$X$, reminiscent of the ``chopping off long
fingers'' construction used by Gromov in the context of surfaces.  As
a consequence, we prove the agreement of the Lusternik-Schnirelmann
and systolic categories of a 2-complex.
\end{abstract}

\maketitle

\tableofcontents

\section{Inequality for unfree complexes}

The homotopy 1-systole, denoted~$\sys(X)$, of a compact metric
space~$X$ is, by definition, the least length of a noncontractible
loop in~$X$.  The first systematic study of this invariant and its
generalisations, by M.~Gromov \cite{Gr1}, built upon earlier
pioneering work of C.~Loewner, P. Pu \cite{Pu}, R. Accola, C.~Blatter,
M.~Berger \cite{Be0}, J.~Hebda~\cite{heb}, and others.  Such work was
recently surveyed in the articles~\cite{CK, KL}.  See \cite{SGT} for
an overview of systolic problems.

The {\em optimal systolic ratio}, denoted~$\SR(X)$, of an
$n$-dimensional complex{}~$X$, is the supremum of the values of the
scale-invariant ratio
\[
\frac{\sys^n(\gmetric)}{\vol_n(\gmetric)},
\]
the supremum being over all piecewise flat metrics~$\gmetric$ on~$X$.

Let~$X$ be a connected finite complex.  Let
\[
f:X \longrightarrow K(\pi_1(X),1)
\]
be its classifying map, \ie a map that induces an isomorphism of
fundamental groups.  Here~$K(\pi,1)$ denotes a connected~$CW$-space
with~$\pi_1(K(\pi,1))=\pi$ and~$\pi_n(K(\pi,1))=0$ for all~$n \geq 2$.
If~$f': X \longrightarrow K'(\pi_1(X),1)$ is another classifying map,
then there exists a homotopy equivalence~$h: K(\pi_1(X),1)\to
K'(\pi_1(X),1)$ such that~$hf\simeq f'$.

\begin{definition}
A complex~$X$ is said to be {\em~$n$-essential\/} if there exists a
classifying map~$X \to K(\pi,1)$ that cannot be homotoped into
the~$(n-1)$-skeleton of~$K(\pi_{1}(X),1)$, \cf \cite[p. 264]{Gr3},
\cite{KR1}.  By the above, {\em every\/} classifying map has this
property.
\end{definition}

In the case~$n=2$, a 2-complex{} is 2-essential if and only if it is
unfree, \cf Theorem~\ref{dim2}.

\begin{definition}
A 2-complex{} is {\em unfree\/} if its fundamental group is not free.
\end{definition}

M.~Gromov~\cite[Appendix~2]{Gr1} proved that every
$n$-essential~$n$-complex~$X$ satisfies a systolic inequality, \ie
there exists a (finite) constant~$C_{n}>0$ such that~$\SR(X) \leq
C_{n}$. Note that a converse is true by~\cite[Lemma~8.4]{Ba}.

When~$X$ is a surface, numerous systolic inequalities are now
available.  They include near-optimal asymptotic upper bounds for the
optimal systolic ratio of surfaces of large genus \cite{Gr1, KS2}, as
well as near-optimal asymptotic {\em lower\/} bounds for large
genus~\cite{BS, KSV}. The genus 2 surface was recently shown to be
Loewner~\cite{KS1}, and moreover admits an optimal systolic inequality
in the CAT(0) class \cite{KS3}. A relative version of Pu's inequality
\cite{Pu} was obtained in~\cite{BCIK1}.  However, many of the existing
techniques, including the entropy technique of~\cite{KS2}, are not
applicable in the more general context of an arbitrary finite
2-complex.

Specifically in dimension 2, M.~Gromov \cite[6.7.A]{Gr1} (note a
misprint in the exponent) showed that every unfree~$2$-dimensional
complex~$X$ satisfies the inequality
\begin{equation}
\label{eq:72}
\SR(X) \leq 10^{4} .
\end{equation}
Contrary to the case of surfaces, where a (better) systolic inequality
can be derived by simple techniques, Gromov's proof of
inequality~\eqref{eq:72} depends on the advanced filling techniques
of~\cite{Gr1}.  Note that the technique of cutting a surface open
along a non-separating loop, and applying the coarea formula to the
distance function from one of the boundary components, does not seem
to generalize to arbitrary~$2$-complexes.

Recently, M.~Gromov \cite{Gr9} remarked that a systolic inequality for
unfree 2-complexes{} should be a consequence merely of the coarea
formula.

We have been able to obtain a {\em uniform\/} inequality for arbitrary
unfree 2-complexes{} using merely the coarea formula, the properties of
the Reeb graph of the distance function, as well as the classical
Grushko decomposition of the fundamental group \cite{St1, SW}.  More
specifically, we construct a kind of a Reeb space ``minimal
model''~$M(X,r)$ for~$X$, reminiscent of the ``chopping off long
fingers'' construction used by M.~Gromov \cite{Gr1} in the context of
surfaces.  We thus improve the constant in~\eqref{eq:72} to the value
12.  Namely, we show the following.

\begin{theorem} 
Every finite unfree~$2$-complex~$X$ satisfies the bound~$\SR(X) \leq
12$.
\end{theorem}

\begin{remark}
Given an arbitrary 2-dimensional complex~$X$, we consider a 2-dimensional 
complex~$Y$ whose fundamental group is the unfree factor of the Grushko
decomposition of the fundamental group of~$X$, \cf
Section~\ref{sec:grush}.  Then there is a map~$f: Y \to X$ inducing a
monomorphism in~$\pi_1$, \cf \cite[Lemma 1.5]{SW}.  Now in the absence
of loose loops (see Section~\ref{three}), for the source space~$Y$, an
inequality~$\SR(Y)\leq 4$ can be established relatively easily by
means of the coarea formula.  Furthermore, one can pushforward
systolic inequalities by such a map~$f$, by a technique pioneered by
I.~Babenko \cite{Ba}, \cf \cite{KR1, KR2}.  This would prove a
systolic inequality for the target space~$X$, as well.  However, the
resulting inequality for~$X$ is not uniform, since pushforward affects
the constant in the inequality.  Namely, the constant is worsened by
the number of faces of~$Y$ in the inverse image of a face of~$X$.  To
overcome this difficulty, we will use an induction on the free index
of Grushko, \cf \eqref{fig}.
\end{remark}

\begin{question}
It is an open question whether all unfree~$2$-complexes{} satisfy
Pu's inequality for~$\rp^{2}$, equivalently if the optimal constant
in~\eqref{eq:72} is~$\frac{\pi}{2}$.
\end{question}

This article is organized as follows.  The corank and the free index
of Grushko~$\FIG(X)$ of a 2-complex{}~$X$ are reviewed in
Section~\ref{sec:grush}.    A useful application of the Seifert-van
Kampen theorem appears in Section~\ref{svk}.  The Reeb graph and its
generalisation called Reeb space are described in Section~\ref{six}.
Section~\ref{sec:tree} introduces a minimal model for~$X$, obtained by
pruning suitable simply connected superlevel sets, as well as
superfluous branches of the Reeb tree.  Technical results on the level
curves of the distance function are presented in Sections~\ref{lc}.
We define loose loops and describe the
intersection of pointed systoles with level curves in
Section~\ref{sec:pointed}.
The dichotomy loose loop/area lower bound is explained in
Section~\ref{three}. In Section~\ref{sec:zero} we prove the bound
$\SR(X)\leq 4$ when~$\FIG(X)=0$, \cf \theoref{theo:balls}.  A suitable
noncontractible level of the distance function is identified in
Section~\ref{sec:fibers}.  In Section~\ref{corank}, using such a
noncontractible level and the coarea formula, we prove a
corank-dependent systolic inequality for unfree~$2$-complexes.  We
compute the Lusternik-Schnirelmann category~$\cat(X)$ of
a~$2$-complex~$X$ in Section~\ref{two}, and show that~$\cat(X)$ and
the systolic category of~$X$ agree.  The proof of the uniform bound
$\SR(X)\leq 12$ occupies Sections \ref{sec:Z} and~\ref{fifteen}.

All complexes{} are assumed to be simplicial, finite, connected, and
piecewise flat, unless explicitly mentioned otherwise.

\section{Corank and free index of Grushko}
\label{sec:grush}

Recall that the {\it corank} of a group~$G$ is defined to be the
maximal rank~$n$ of a free group~$F_n$ admitting an epimorphism~$G \to
F_n$ from~$G$.  Clearly, every finitely generated group has finite
corank, and therefore each finite~$CW$-space has fundamental group of
finite corank.

Grushko's theorem~\cite{St1, SW} asserts that every finitely generated
group~$G$ has a decomposition as a free product of subgroups
\begin{equation}
\label{41b}
G=F_{p}*H_{1}*\dots*H_{q}
\end{equation}
such that~$F_{p}$ is free of rank~$p$, while every~$H_{i}$ is
nontrivial, non isomorphic to~$\Z$ and freely indecomposable.
Furthermore, given another decomposition of this sort, say
$G=F_{r}*H'_{1}*\dots*H'_{s}$, one necessarily has~$r=p$,~$s=q$ and,
after reordering,~$H'_{i}$ is conjugate to~$H_{i}$.

\begin{definition}
\label{41y}
We will refer to the number~$p$ in decomposition \eqref{41b} as the
{\em free index of Grushko\/} of~$G$, denoted~$\FIG(G)$.
\end{definition}

Thus, every finitely generated group~$G$ with~$\FIG(G)=p$ can be
decomposed as
\begin{equation}
G = F_{p} * H_{G},
\end{equation}
where~$F_{p}$ is free of rank~$p$ and~$\FIG(H_G)=0$.  The
subgroup~$H_{G}$ is unique up to isomorphism.  Its isomorphism class
is called the unfree factor of (the isomorphism class of)~$G$.  If~$X$
is a finite complex, we set
\begin{equation}
\label{fig}
\FIG(X)=\FIG(\pi_1(X))
\end{equation}
by definition.  Note that every finitely generated group~$G$ satisfies
$\FIG(G)\leq \corank(G)$.

\begin{definition} 
\label{def:split}
Let~$G$ be a finitely generated group. We say that an element~$g\in G$
{\em splits off\/} if~$G$ admits a decomposition as a free product
with the infinite cyclic group generated by~$g$.
\end{definition}

\section{An application of the Seifert-van Kampen theorem}
\label{svk}

We need the following well-known fact about~$CW$ complexes, \cf
\cite[Prop. I.3.26]{Rud}.

\begin{prop}\label{cofib}
Let~$(X,A)$ be a~$CW$-pair. If~$A$ is a contractible
space, then the quotient map~$X \to X/A$ is a homotopy
equivalence.
\qed
\end{prop}

\begin{cory}
\label{cone}
Consider a~$CW$-pair~$(X,A)$ and a~$CW$-space~$X\cup CA$ where~$CA$ is
the cone over~$A$.  Then the quotient map
\[
X\cup CA\to (X\cup CA)/CA=X/A
\]
is a homotopy equivalence.  In particular,~$X\cup CA$ and~$X/A$ are
homotopy equivalent,~$X\cup CA \simeq X/A$.  \qed
\end{cory}

\begin{lemma}\label{pi-exc}
Let~$(X,A)$ be a~$CW$-pair with~$X$ and~$A$ connected.  Then the
quotient map~$q: X \to X/A$ induces an epimorphism of fundamental
groups. Furthermore, if the inclusion~$j: A \to X$ induces the zero
homomorphism~$j_*: \pi_1(A) \to \pi_1(X)$ of fundamental groups, then
the quotient map~$q: X \to X/A$ induces an isomorphism of fundamental
groups.
\end{lemma}

\begin{proof}
By the Seifert--van Kampen theorem, the inclusion \[ i:X \subset X\cup
CA
\]
induces an epimorphism~$i_*: \pi_1(X) \to \pi_1(X \cup CA)$ of
fundamental groups, and~$i_*$ is an isomorphism if~$j_*$ is the zero
map.  Finally, by Corollary~\ref{cone}, the quotient map~$X \cup CA \to X/A$
is a homotopy equivalence, while the map~$q$ is the composition~$X
\subset X \cup CA \to X/A$, proving the lemma.
\end{proof}

\begin{lemma}\label{wedge}
Let~$A=\{a_0, a_1, \ldots, a_k \}$ be a finite subset of a
connected~$CW$-space~$X$.  Then~$X \cup CA$ is homotopy equivalent to
the wedge of~$X$ and~$k$ circles,~$X\cup CA \simeq X \vee S^1_1\vee
\cdots \vee S^1_k$. In particular,
$$
\pi_1(X\cup CA)=\pi_1(X)*F_k
$$ 
where~$F_k$ is the group of rank~$k$. In other words, the free
index of Grushko of~$\pi_1(X\cup CA)$ is at least~$k$.
\forget
Furthermore, the loop homotopy class of the loop~$\overline{\alpha}$, the image 
under the
quotient map~$X \longrightarrow X/A$, of an embedded arc~$\alpha$ in
$X$ joining two points of~$A$, splits off in~$\pi_1(X/A)$.
\forgotten
\end{lemma}

\begin{proof}
Let~$X^{(1)}$ be the 1-skeleton of~$X$. Without loss of generality, we
can assume (subdividing~$X$) that~$A\subset X^{(1)},\, \alpha\subset
X^{(1)}$ and that~$X^{(1)}$ is connected. Since~$X^{(1)}$ is connected,
we can find a tree~$T\subset X^{(1)}$ that contains~$\alpha$ and~$A$.
Since~$T$ is contractible, we have
\begin{equation*}
\begin{aligned}
X\cup CA &\simeq (X\cup CA)/T\simeq (X/T) \vee (CA/A)\\
& \simeq X\vee (CA/A) \simeq X \vee
S^1_1\vee \cdots \vee S^1_k
\end{aligned}
\end{equation*}
from \propref{cofib}.
\forget
The loop~$\overline{\alpha}$ can be represented by a loop~$c = \alpha
\cup \beta \subset X \cup CA$ under the homotopy equivalence~$X \cup
CA \simeq X/A$ of Corollary~\ref{cone}, where~$\beta$ agrees with the
unique embedded arc of~$CA$ with the same endpoints as~$\alpha$.
Since~$\alpha \subset T$ and its endpoints are distinct, the loop~$c$
can be represented by a noncontractible simple loop of~$CA/A$ under
the homotopy equivalence~$X \cup CA \simeq X \vee (CA/A)$.  Thus, the loop 
homotopy 
class of the 
loop~$c$ (and so~$\overline{\alpha}$) splits off.
\forgotten
\end{proof}

Let~$X$ be a finite connected complex{} and let~$f: X \to \R^+$ be a
function on~$X$.  Let
\[
[f\le r]:=\{x\in X\bigm| f(x)\le r\} \quad \text{ and } \quad [f\ge
r]:=\{x\in X\bigm| f(x)\ge r\}
\]
denote the sublevel and superlevel sets of~$f$, respectively.

\begin{definition} \label{def:coalesce}
Given~$r\in \R^+$, suppose that a single path-connected component
of the superlevel set~$[f\ge r]$ contains
$k$ path-connected components of the level set~$f^{-1}(r)$.
Then we will say that the~$k$ path-connected components {\em
coalesce forward}.
\end{definition}

For future needs, recall that the connected components of any complex{}
are path-connected.

\begin{lemma} 
\label{42}
Given~$r>0$, assume that the pair~$\left([f\ge r], f^{-1} _{\phantom
{I}} (r)\right)$ is homeomorphic to a~$CW$-pair. Suppose that the
set~$[f\le r]$ is connected and that~$k+1$ connected components
of~$f^{-1}(r)$ coalesce forward. Then the corank of~$\pi_1(X)$ is at
least~$k$.  Furthermore, if the inclusion
\[
[f\le r]\subset X
\]
of the sublevel set~$[f \leq r]$ induces the zero homomorphism of
fundamental groups, then~$\FIG(X)\geq k$.
\end{lemma}

\begin{proof} Let~$C_0, C_1 \ldots, C_k$ be the distinct
components of~$f^{-1}(r)$ that coalesce forward, and let~$Z$ be the
component of~$[f \ge r]$ that contains all~$C_i$'s.  Then~$[f \ge
r]=Z\sqcup W$ where~$W$ is the (finite) union of all components of~$[f
\ge r]$ other than~$Z$.  Let~$Y=Z/\sim$ where~$x\sim y$ if and only
if~$x,y$ belong to the same component~$C_i$ of~$f^{-1}(r)$. Let~$a_i$
be the image of~$C_i$ under the quotient map~$Z \to Y$ and
let
\[
A=\{a_0, \ldots a_k\}\subset Y.
\]
We set
\[
W'=W/(W\cap f^{-1}(r)).
\]
Then we have
\[
X/[f\le r]= W'\vee Z/(\cup C_i)=W'\vee Y/A\simeq W'\vee Y\vee
S^1_1\vee \cdots \vee S^1_k ,
\]
where the last equivalence follows from Corollary~\ref{cone} and
\lemref{wedge}.  By the Seifert--van Kampen Theorem, we have
$\FIG(X/[f\le r]) \geq k$. In particular, there exists an epimorphism
\[
\pi_1(X/[f\le r])\to F_k.
\]
By \lemref{pi-exc}, the quotient map~$X \to X/[f\le r]$ induces an
epimorphism of fundamental groups.  Hence~$\corank(\pi_1(X))\geq
k$. Finally, if the inclusion~$[f\le r]\subset X$ induces the zero
homomorphism of fundamental groups, then~$\pi_1(X/[f\le r])$ is
isomorphic to~$\pi_1(X)$ by \lemref{pi-exc}, proving the lemma.
\end{proof}

\section{Reeb graph of a piecewise flat complex}
\label{six}

\begin{definition}
Let~$X$ be a finite connected complex.  Consider a function~$f:X
\longrightarrow \R^+$.  Let~$r_0>0$ be a real number.  The Reeb space
of~$f$ in range~$r_0$, denoted~$\Reeb(f,r_0)$, is defined as the
quotient
\[
\Reeb(f,r_0) = X/\sim
\]
where we have equivalence~$x \sim y$ if and only if~$f(x)=f(y)\leq
r_0$, and~$x$ and~$y$ lie in the same connected component of the level
set~$f^{-1}(f(x))$, \cf \cite{Re}.  The Reeb space~$\Reeb(f,\infty)$
in ``full range'' will be denoted simply~$\Reeb(f)$.
\end{definition}

The following fact is a consequence of standard results on the
triangulation of semialgebraic functions 
\cite[\S 9]{bcr}, \cite{Sh}.

\begin{proposition}
\label{prop:reeb}
Let~$X$ be a finite,~$2$-dimensional, piecewise flat complex.  Then
the Reeb space~$\Reeb(f)$ of the distance function~$f=f_p$ from a
point~$p\in X$ is a finite graph.

Furthermore, the finite graph~$\Reeb(f)$ can be subdivided so that the
natural map~$X \to \Reeb(f)$ yields a trivial bundle over the interior
of each edge of~$\Reeb(f)$.
\end{proposition}

The Reeb graph was used in \cite{GL, Ka1b} to study 3-manifolds of
positive scalar curvature.  Other applications are discussed in
\cite{CP}.

\begin{remark}
A more precise description (than that predicted by semireal algebraic
geometry) can be given of the level curves themselves, \cf
Section~\ref{lc}.
\end{remark}

From now on,~$X$ and~$f$ will be as in Proposition~\ref{prop:reeb}.
The Reeb space~$\Reeb(f,r)$ in range~$r$ can be thought of as a
``hybrid'' space, consisting of two pieces.  One piece is the Reeb
graph~$\Reeb(f|_B)$ of the ball~$B=B(p,r)$ of radius~$r$ centered
at~$p$.  The other piece~$[f \geq r]$ is the complement, in~$X$, of
the open ball of radius~$r$, attached to the graph by a map~$\mu$
which collapses each connected component of the level
curve~$f^{-1}(r)$ to a point:
\begin{equation}
\label{hybrid}
\Reeb(f,r) = \Reeb \left( f|_{B(p,r)} \right) \cup_\mu [f \geq r].
\end{equation}

The Reeb graph~$\Reeb(f)$ of~$f$ is endowed with the length structure
induced from~$X$.  Let~$\overline{p}\in \Reeb(f)$ be the image of~$p$.
The ball
\[
B(p,r) \subset X
\]
of radius~$r < \tfrac{1}{2}\sys(X)$ projects to the
ball~$B(\overline{p},r) \subset \Reeb(f)$.

\begin{lemma}
\label{25b}
Let~$r < \tfrac{1}{2}\sys(X)$.  The subgraph~$B(\overline{p},r)\subset
\Reeb(f)$ is a tree, denoted~$T_r$.  Thus, the decomposition
\eqref{hybrid} becomes
\begin{equation}
\Reeb(f,r) = T_r \cup_\mu [f \geq r].
\end{equation}
\end{lemma}

\begin{proof}
Every edge of~$\Reeb(f)$ lifts to a minimizing path in~$X$, given by a
segment of the minimizing path from~$p$.  Given an embedded
loop~$\gamma \subset \Reeb(f)$, we can lift each of its edges to~$X$.
The endpoints of adjacent edges lift to a pair of points lying in a
common connected component of a level curve of~$f$ (by definition of
the Reeb graph), and can therefore be connected in~$X$ by a path which
projects to a constant path in~$\Reeb(f)$.  We thus obtain a loop
in~$X$ whose image in~$\Reeb(f)$ is homotopic to~$\gamma$.  The lemma
now follows from the fact that the inclusion~$B(p,r) \subset X$
induces the trivial homomorphism of fundamental groups, \cf
Remark~\ref{ball2}.
\end{proof}

\forget Suppose that the graph~$B(\overline{p},R)$ has a nontrivial
cycle.  Then, we can construct a loop~$\gamma \subset B(p,R)$ formed
of two distance minimizing arcs joining~$p$ to some point~$q$ that
projects to a noncontractible loop of~$B(\overline{p},R)$.  The
loop~$\gamma$ is contractible in~$X$.  Therefore, there exists a
simplicial map~$s:D \longrightarrow X$ from the unit disk~$D$ that
takes~$\partial D$ to~$\gamma$.  Let~$x$ and~$y$ be two points
of~$\gamma$ that lie in two different components of~$f^{-1}(r)$.  We
identify~$\gamma$ to~$\partial D$.  Thus,~$x$ and~$y$ lie in two
different components of~$D_{r} = (f \circ s)^{-1}(r)$.  In particular,
there exists an arc of~$D$ from~$p$ to~$q$ that does not
intersect~$D_{r}$.  Therefore,~$q$ is closer to~$p$ than~$x$ or~$y$.
Hence a contradiction.  \forgotten

\begin{lemma}
\label{24}
Let~$r < \tfrac{1}{2}\sys(X)$.  
Then the natural projection map~$h: X\to
\Reeb(f,r)$ induces an isomorphism of fundamental groups,
\begin{equation}
h_*: \pi_1(X) \to \pi_1(\Reeb(f,r)).
\end{equation}
\end{lemma}

\begin{proof}
The obvious map~$g: X/B(p,r) \to \Reeb(f,r)/T_r$ (induced by~$h$) is a 
homeomorphism. Consider the commutative diagram
\[
\CD \pi_1(X) @>h_*>> \pi_1(\Reeb(f,r))\\ @VVV @VVV\\ \pi_1(X/B(p,r))
@>g_*>> \pi_1(\Reeb(f,r)/T_r) \endCD
\]
The vertical maps are isomorphisms by \lemref{pi-exc} and
Remark~\ref{ball2}, while~$g_*$ is induced by a homeomorphism.
\end{proof}

\forget
\begin{lemma}\label{25}
Let~$\gamma$ be a pointed systolic loop at~$p \in X$, and let
$L=\length(\gamma)=\sys(X,p)$. Let~$r$ be a real number satisfying 
\begin{equation}
\eqref{rr}
L - \sys(X) <2r< L.
\end{equation}
If~$\gamma$ meets two distinct connected components of the level curve
$S(r)\subset X$, then~$\gamma$ is loose.
\end{lemma}

\begin{proof}
Consider the natural projection~$X \to \Reeb(f,r)$. Then, by
\lemref{wedge} and \lemref{lem:simple}, the image of~$\gamma$ in
$\Reeb(f,r)$ is a loose loop (the lower bound in \eqref{rr} is needed
to guarantee imbeddedness, see Lemma~\ref{lem:simple}(ii) below).  On
the other hand, the natural projection induces an isomorphism of
fundamental groups by \lemref{24}.
\end{proof} 
\forgotten

\section{A minimal model and tree energy}
\label{sec:tree}

We describe a pruning of the Reeb space which results in a kind of a
``minimal model'', \cf \eqref{mm}, for our 2-complex{}~$X$.

Let~$r< \tfrac{1}{2}\sys(X)$.  Consider the function~$\bar f:
\Reeb(f,r) \to \R$ naturally defined by~$f$.  Denote
by~$\mathcal{C}\subset \Reeb(f,r)$ the union of the connected components
$C$ of the subset~$[\bar f \geq r]$ with the following two properties:
\begin{enumerate}
\item
$C$ is simply connected;
\item
$C$ is attached to the tree~$T_{r} \subset \Reeb(f,r)$ at a single
leaf (\ie vertex attached to precisely one edge).
\end{enumerate}

The pullback of the standard 1-form~$dr$ on~$\R$ by the map
\[
\Reeb(f)\to \R
\]
back to an edge of~$\Reeb(f)$ defines a positive direction on each
edge.  We prune the tree~$T_r$ by removing edges of the following two
types:

\begin{enumerate}
\item
edges that do not reach the set~$[\bar f \geq r]$ when we follow the
positive direction;
\item
edges that lead only to components in~$\mathcal C \subset
[\bar f \geq r]$.  
\end{enumerate}

Denote by~$T'_r \subset T_r$ the resulting pruned tree.  In other
words,~$T'_r$ is the union of embedded paths connecting~$p$ to
components of the complement~$[\bar f \geq r] \setminus \mathcal C$.

A vertex of~$T'_r$ is by definition a point whose open neighborhood
in~$T'_r$ is not homeomorphic to an open interval.  An edge is a
connected component of the complement of the set of vertices.

Define the {\em energy\/}~$E(\Gamma)$ of a graph~$\Gamma$ to be the
sum of the squares of the lengths of its edges.  The height of a tree
with a distinguished root leaf~$p$ is the {\em least\/} distance
from~$p$ to any other leaf of the tree.

\begin{proposition}
\label{21}
The energy of a tree~$\Gamma$ of height~$h$ satisfies the
bound~$E(\Gamma)\geq \tfrac{1}{2}h^2$.
\end{proposition}

\begin{proof}
For a simplest tree shaped as the letter Y, the least energy is given
by the metric for which the bottom interval is twice as long as each
of the top two intervals.  Arguing inductively, we see that the lower
bound is attained by the infinite tree where the length of an edge is
halved after each branching.
\end{proof}

\begin{definition}
We introduce a kind of a minimal model~$M(X,r) \subset \Reeb(f,r)$
for~$X$ by setting
\begin{equation}
\label{mm}
M(X,r) = T'_r \cup_\mu ([f\geq r] \setminus {\mathcal C}).
\end{equation}
\end{definition}

\begin{proposition}
\label{33}
The map~$X \to M(X,r)$ collapsing all superfluous material to the
vertex of~$T_r$ where it is attached to the pruned tree~$T'_r \subset
T_r$, induces an isomorphism~$\pi_1(X) \simeq \pi_1(M(X,r))$.
\end{proposition}

\forget
\begin{prop}
\label{inner}
If~$X$ admits a systolic loop~$\gamma$ at~$p$ that intersects two
different components of~$S(r)$ with~$\sys(X,p)-\sys(X)<2r<\sys(X,p)$,
then~$M(X,r)$ admits a non-root edge.
\end{prop}

\begin{proof} 
Let~$a, b$ be two different points where~$\gamma$ meets~$S(r)$.  Then
$M(X,r)$ has two different edges, and at least one of them is non-root.
\end{proof}
\forgotten

\begin{definition}
\label{df:inner}
If the pruned tree~$T'_r$ has exactly one edge attached to its root leaf, we 
call this edge the {\em root edge} of~$T'_r$.
If~$T'_r$ has more than one edge attached to its root leaf, the root edge is 
defined as the trivial edge reduced to the root leaf of~$T'_r$.
\end{definition}

\begin{proposition}
\label{32}
Let~$X$ be an unfree~$2$-complex.  Let~$e \subset T'_r$ be an open non-root
edge. Consider the space
\[
Y= M(X,r) \setminus e
\]
obtained from the minimal model by removing~$e$.  Then there exists an
unfree connected component~$Y_*\subset Y$ satisfying
\begin{equation}
\label{52f}
\FIG(Y_*) \leq \FIG(X)-1.
\end{equation}
\end{proposition}

\begin{proof}
If removing~$e$ does not disconnect the space~$\Reeb(f,r)$, then, by
\lemref{wedge}, we have
\[
\pi_{1}(M(X,r)) \simeq \pi_{1}(\Reeb(f,r)) \simeq \pi_{1}(Y)*\Z,
\]
proving equality in~\eqref{52f}.  If removing~$e$
disconnects~$\Reeb(f,r)$, then~$Y$ decomposes as~$Y_{-} \cup Y_{+}$,
where each of~$Y_{\pm}$ is connected.  Neither component is simply
connected by definition of the pruned tree.  We have~$\pi_{1}(X)
\simeq \pi_{1}(Y_{-}) * \pi_{1}(Y_{+})$ by the Seifert-van Kampen
theorem.  Hence~$\FIG(X)=\FIG(Y_-) + \FIG(Y_+)$.  If both~$Y_+$
and~$Y_-$ are unfree, we choose the one with the least~$\FIG$.  If one
of them is free, then the other necessarily satisfies \eqref{52f}.
\end{proof}

\section{Level curves} 
\label{lc}

In this section, we present technical results concerning the level
curves of the distance function on a piecewise flat 2-complex.  The
2-simplices of such a complex{} will be referred to as {\em faces}.

\begin{theorem}
\label{46lc}
A level curve of the distance function from a point in a finite
piecewise flat~$2$-dimensional complex~$X$ is a finite union of
circular arcs and isolated points.
\end{theorem}

The claim is obvious for level curves at sufficiently small distance
from a vertex~$p\in X$.  Such a level curve is known as the {\em
link\/} of~$X$ at~$p$, \cf \cite[\S 2]{Mu}.

The structure of the proof is 3-fold:
\begin{enumerate}
\item
a finiteness result on the number of geodesic arcs of bounded length
joining a pair of points of~$X$;
\item
the construction of a finite graph~$\Lambda_p$, \cf \eqref{51lc},
containing a level curve~$S$;
\item
an argument showing that~$S$ is a subgraph of~$\Lambda_p$.
\end{enumerate}

\begin{lemma}
\label{52}
Let~$C>0$.  Then there are finitely many geodesic arcs of length at
most~$C$ between any pair of points of a finite piecewise flat
$2$-complex~$X$.
\end{lemma}

\begin{proof}
Let~$X'$ be the complement of the set of vertices in~$X$.  Let~$p,q\in
X'$.  By the flatness of the complex, there is at most one geodesic
path from~$p$ to~$q$ in a given homotopy class of such paths in~$X'$.
Note that~$X'$ is homotopy equivalent to a compact 2-complex, by
removing an open disk around each vertex.  The lemma now follows from
the finiteness of the number of homotopy classes of bounded length in
a compact complex.
\end{proof}

\begin{remark}
\label{53}
Consider a compact geodesic segment~$\gamma\subset \widetilde X'$ in
the universal cover~$\widetilde X'$ of~$X'$.  Then the union of all
faces meeting~$\gamma$, includes a flat strip containing~$\gamma$.
\end{remark}

Continuing with the proof of Theorem~\ref{46lc}, note that the link at
$q$ is a graph which can be identified with the set of unit tangent
vectors at~$q$.  To describe the graph~$\Lambda$ of step 2, we refine
the link at~$q$ as follows.  Consider the tangent vectors of the
geodesics of length smaller than the diameter of~$X$, connecting~$q$
to all vertices of~$X$.  We subdivide the link recursively.  At each
step, we add to the link the vertices corresponding to such tangent
vectors, so that the edge containing such a vertex is split into two.
By Lemma~\ref{52}, the resulting graph~$V_q$ is finite.

Choose a number~$\rho$ satisfying~$0<\rho \leq \diam(X)$.  Given an
edge~$e$ of the refined link~$V_{q}$ at~$q$, consider the
corresponding pencil of geodesics, of length~$\rho$, issuing
from~$q$. By Remark~\ref{53}, their endpoints trace out a single
compact (possibly self-intersecting) circular arc, denoted
\begin{equation*}
\CA(q,e,\rho) \subset X,
\end{equation*}
of ``radius''~$\rho$ and geodesic curvature~$\tfrac{1}{\rho}$, at
least if~$X$ is a surface.  In general, the structure
of~$\CA(q,e,\rho)$ may be more complicated due to branching at a point
where a geodesic issuing from~$q$ encounters a 1-cell of~$X$, but in
any case the number of circular arcs forming the graph~$\CA$ can be
controlled in terms of the total number of faces.  We set~$\CA (q, e,
\rho)= \emptyset$ if~$\rho<0$.

Now let~$p\in X$ be a vertex, and choose a number~$r$ satisfying~$0< r
\leq \diam(X)$.  Let~$(q_i)$ be an enumeration of the vertices
of~$X$, and~$(e_{ij})$ an enumeration of the edges of the refined
link~$V_{q_i}$.

\begin{lemma}
\label{54}
Let~$S=S(r)\subset X$ be the level~$r$ curve of the distance function
from~$p\in X$, namely,~$S=\{x \in X \bigm |\dist (p,x)=r\}$.  Consider
the finite union
\begin{equation}
\label{51lc}
\Lambda_p =\Lambda_p(r)= \bigcup_{i,j} \CA \left(q_i,
e_{ij}^{\phantom{I}}, r-\dist(p,q_i) \right).
\end{equation}
Then~$\Lambda_p$ is a finite graph, containing~$S$, and at distance at
most~$r$ from~$p$.
\end{lemma}

\begin{proof}
To describe the graph structure of~$\Lambda_p$, note that a pair of
circular arcs in a given simplex of~$X$ meet in at most a pair of
points.  This elementary algebraic-geometric observation is thus the
basis of the proof of Theorem~\ref{46lc}; \cf Remark~\ref{59}.

Whether the arcs are transverse or tangent, we include the common
points as vertices of~$\Lambda_p$. Also, if a circular arc of
$\Lambda_p$ meets an edge of~$X$, then the common point is declared to
be a vertex of~$\Lambda_p$.  \forget Finally,~$\Lambda_p$ may contain
an entire circle.  In such case we add some vertices to the circle so
as to turn it into a subgraph.  \forgotten
\end{proof}

To complete the proof of Theorem~\ref{46lc}, we need a notion of a cut
locus.  The notion of a cut locus for a complex{} is ill-defined.
However, we define a ``local'' analogue~$\CL_p$ of the cut locus from
a point~$p\in X$, where~$\CL_p$ is inside an open face~$\Delta\subset
X$, as the set of points~$q\in \Delta$ with at least a pair of
distinct unit tangent vectors at~$q$ to minimizing geodesics from~$p$.

\begin{remark}
\label{55}
The cut locus on a surface is generally known to be a graph, and the
same can be shown for our~$\CL_p$.  The vertices of~$\CL_p$ seem to be
related to the critical points of the distance function from~$p$ in
the sense of Grove-Shiohama-Gromov-Cheeger \cite{Che}.  We will not
pursue this direction.
\end{remark}

\begin{lemma} 
\label{lem:relclosed}
The set~$\CL_{p}\subset X$ is relatively closed inside each open face
of~$X$.
\end{lemma}

\begin{proof}
Consider a sequence of points~$\{x_n\}$ in~$\CL_{p}$ converging
to~$x\in \Delta$.  Let~$u_n$ be a unit vector at~$x_n$, and tangent to
a minimizing geodesic from~$p$.  By choosing a subsequence, we can
assume that the minimizing geodesics corresponding to the
vectors~$u_n$ lead to the same vertex~$A\in X$, and, furthermore, are
homotopic as relative paths from~$A$ to~$\Delta$ for the pair
$(\{A\}\cup X', \Delta)$, where~$X'$ is the complement of the set of
vertices in~$X$.  Now, we use a kind of a developing map and consider
the face~$\Delta$ as lying in the Euclidean plane~$\R^2$.  Then there
is a point~$\tilde A\in \R^2$ such that the circular arcs of
$\Lambda_p$ corresponding to the vectors~$u_n$ at~$x_n\in \Delta
\subset \R^2$ are arcs of concentric circles with common
center~$\tilde A$.

If~$\{v_n\}$ is another sequence of tangent vectors at~$\{x_n\}$, we
similarly obtain points~$B\in X$ and~$\tilde B\in \R^2$ for the
sequence~$v_n$.  Now suppose the two sequences have a common limit
vector~$\lim_n u_n=\lim_n v_n$.  Then the points~$\tilde A, \tilde B,
x \in \R^2$ are collinear, while~$x\not\in [\tilde A,\tilde B]$.
Hence
\[
\begin{aligned}
d(\tilde A, \tilde B) & = | d(\tilde A, x) - d(\tilde B, x)| \\ & = |
d(A, p) - d(B, p)|.
\end{aligned}
\]
Since the geodesics leading to~$x_n$ are assumed minimizing, we have
\[
\begin{aligned}
|d(\tilde A,x_{n}) - d(\tilde B,x_{n})| & = | d(A, p) - d(B, p)| \\
& = d(\tilde A, \tilde B).
\end{aligned}
\]
Therefore, the points~$\tilde A, \tilde B, x_n \in \R^2$ are
collinear, as well, and~$x_n \not\in [\tilde A, \tilde B]$.  Thus
$u_n=v_n$, proving the lemma.
\end{proof}

The proof of Theorem~\ref{46lc} is now completed by means of the
following lemma.

\begin{lemma} \label{lem:subgraph}
The level curve~$S$ is a subgraph of the finite graph~$\Lambda_p$.
\end{lemma}

\begin{proof}
Since~$S$ meets the 1-skeleton in a finite number of points, it
suffices to examine the behavior of~$S$ inside a single face~$\Delta$.

At least two circular arcs corresponding to edges of~$\Lambda_p$ pass
through each point of~$\CL_{p} \cap \Lambda_p$.  Therefore, the points
of~$\CL_{p} \cap \Lambda_p$ are vertices of the graph~$\Lambda_p$.  In
particular, every open edge~$e \subset \Lambda_p$ is disjoint
from~$\CL_{p}$.

Every point~$x\in e \cap S$ is away from the relatively closed
set~$\CL_{p}$, \cf Lemma~\ref{lem:relclosed}.  Hence the level curves
(of the distance function from~$p$) in the neighborhood of~$x$ form a
family of concentric circular arcs.  The circular arc of this family
passing through~$x$ is clearly contained in~$S$, but also in~$e$.
Otherwise, it would intersect transversely the circular arc~$e$ and
the point~$x$ would lie in~$\CL_{p}$.  We conclude that~$e \cap S$ is
open in~$e$.

On the other hand, the level curve~$S$ of a distance function is a
closed set, hence the intersection~$e \cap S$ is closed in~$e$.
Therefore, this intersection coincides with~$e$ or is empty.  It
follows that~$S\subset \Lambda_p$ is a subgraph.
\end{proof}
 
\forget
In other words,~$\Lambda \setminus S$ is a disjoint union of open edges and 
vertices of~$\Lambda$. In future we call these edges of~$\Lambda$ {\em 
excessed}. 
\forgotten

Given a point~$p\in X$, set~$f(x)=\dist (p,x)$.

\begin{corollary}\label{homeo}
The triangulation of~$X$ can be refined in such a way that the sets
$[f \leq r]$,~$f^{-1}(r)$, and~$[f \geq r]$ become~$CW$-subspaces
of~$X$.
\end{corollary}

\begin{proof}
We add the graph~$S=f^{-1}(r)$ to the 1-skeleton of~$X$.  Some of the
resulting faces may not be triangles, therefore a further (obvious)
refinement may be necessary.  The sets~$[f \leq r]$ and~$[f \geq r]$
are connected components of the complement~$X\setminus S$, and hence
$CW$-subspaces.
\end{proof}

\forget
\begin{lemma}\label{meet}
Let~$\Delta$ be an open face in~$X$. Suppose that there are~$2$
circular arcs~$C_1, C_2$ such they either intersect transversally at
$x\in \Delta$ or are tangent at~$x\in \Delta$ with the same direction
of concavity. Let~$e_i, i=1,2,3,4$ be the edges of~$\Lambda$ contained
in~$C_i, i=1,2$ and having~$x$ as an end-point. Then at least two of
these edges are excessed. Furthermore, if~$\{C_1, \ldots C_r\}$ is the
set of all circular arcs passing through the vertex~$x$ of
$\Lambda\cap \Delta$, then at least two of the edges ending at~$x$ are
not excessed.
\end{lemma}

\begin{proof} This first claim follows directly from obvious concavity 
arguments. To prove the second claim, first notice that if~$x\in S(r),
r>0$then at least one edge ending at~$x$ must be in~$S(r)$. Indeed,
$x$ appears as the intersection of circular arcs, and so at least one
of these arcs must have arbitrary close to~$x$ points in~$S(r)$. Now,
suppose that there is exactly one edge~$e$ that ends at~$x$. Then, for
an arbitrary small neighborhood~$U$ of~$x$, we can find to points
$y,z\in U$ such that~$f(y)<r <f(z)$. We can assume that~$U\setminus
\ov e$ is connected, and therefore there exists a point~$w\in
U\setminus \ov e$ with~$f(w)=r$. This implies the existence of one
more edge in~$S$ ended in~$x$.
\end{proof}

\begin{cory} 
If~$x\in S\cap \Delta$ is not a vertex of~$S$ where two tangent
circular arcs have opposite directions of concavity, then~$x$ belongs
at most two edges of~$S$.  
\qed
\end{cory}

Let~$\mathcal{C}_{r}$ be the collection of the edges (circular
arcs) of the finite graph~$\Lambda$ defined in~\eqref{51lc}.

We say that two circular arcs in the plane are {\it independent} if
they do not lie in the same circle. We also say that two edges
of~$\Lambda$ are {\it independent} if they lie in the same~$2$-simplex
and are independent as circular arcs.  
\forgotten

\begin{remark}
\label{59}
Alternatively, Corollary~\ref{homeo} (but {\em not\/}
Theorem~\ref{46lc}) can be deduced from standard results in real
semialgebraic geometry, as follows, \cf \cite{bcr}.  First, note
that~$X$ can be embedded into some~$\R^{N}$ as a semialgebraic set and
that the distance function~$f$ is a semialgebraic function on~$X$.
Thus, the level curve~$f^{-1}(r)$ is a semialgebraic subset of~$X$
and, therefore, a finite graph, \cf proof of Lemma~\ref{54}.
\end{remark}

\section{Loose loops and pointed systoles}
\label{sec:pointed}

\forget
We say that~$\gamma$ is a {\em loose loop\/} if one of the following
two equivalent conditions is satisfied:
\begin{enumerate}
\item
the loop homotopy class~$[\gamma]$ of~$\gamma$ splits off in
$\pi_1(X,p)$;
\item
there is a 2-complex~$Y$ and a map~$X\to Y \vee S^1$ inducing an
isomorphism~$\varphi :\pi_1(X,p) \to \pi_1(Y, y_0) * \Z$ where the
image of~$[\gamma]$ is a generator of the free factor:~$\varphi
([\gamma]) = 1\in \pi_1(S^1)=\Z$.
\end{enumerate}
\forgotten

\forget
\begin{lemma}
\label{25}
Let~$\gamma \subset \Reeb(f,r)$ be an embedded loop in the Reeb space
in range~$r$, and assume~$\gamma$ decomposes as the union~$\gamma
\gamma' \cup \gamma''$ where the paths~$\gamma' \subset \Reeb(f|_B)$
and~$\gamma'' \subset [f \geq r]$ have common endpoints.
Then~$\gamma$ is loose.
\end{lemma}

\begin{proof}
The path~$\gamma''\subset \Reeb(f,r)$ is contractible.  Hence we can
think of the fundamental group as the relative
group~$\pi_1(\Reeb(f,r), \gamma'')$.
\end{proof}
\forgotten

\begin{definition}
\label{ball}
Let~$p\in X$.  A shortest noncontractible loop of~$X$ based at~$p$ is
called a {\it pointed systolic loop\/} at~$p$. Its length, denoted
by~$\sys(X,p)$, is called the {\it pointed systole\/} at~$p$.
\end{definition}

\begin{remark}
\label{ball2}
Alternatively,~$\sys(X,p)$ could be defined as twice the upper bound
of the reals~$r > 0$ such that induced map~$\pi_1(B(p,r))\to \pi_1(X)$
is zero. In other words, every loop in~$B(p,r)$ is contractible
in~$X$.
\end{remark}

\forget
\begin{prop}\label{ball}
Let~$B(p,r)$ be the ball centered at~$p \in X$ of radius~$r <
\frac{1}{2} \sys(X,p)$.  Then the inclusion~$B(p, r)\subset X$ induces
the zero homomorphism of fundamental groups.
\end{prop}

\begin{proof}
Suppose the contrary and consider all loops at~$p$ that are
noncontractible in~$X$. Let~$L$ denote the minimum of length of such
loops; then~$L\ge \sys(X,p)$.  Let~$\gamma$ be the noncontractible
loop in~$(X,p)$ that realizes this minimum length. Let~$a\in B(p,r)$
be the point of~$\gamma$ that divides~$\gamma$ into two arcs
$\gamma_1$ and~$\gamma_2$ of the same length~$L/2$. Consider the
smallest geodesic path~$c$ that joins~$p$ to~$a$. It has the length
$d<r$. Since at least one of the curves~$\gamma_1\cup c_-$ or~$c\cup
\gamma_2$ is noncontractible, we conclude that~$d+L/2\ge L$, \ie~$d\ge
L/2$ (here~$c_-$ denotes the path~$c$ with the opposite
orientation). Thus
$$
\sys(X,p)>2r\ge 2d\ge L \ge\sys(X,p).
$$
That is a contradiction.
\end{proof} 
\forgotten

The following lemma describes the structure of a pointed systolic loop.

\begin{lemma} 
\label{lem:simple}
Let~$\gamma$ be a pointed systolic loop at~$p \in X$, and let
$L=\length(\gamma)=\sys(X,p)$.
\begin{enumerate}
\item[(i)] The loop~$\gamma$ is formed of two distance-minimizing
arcs, starting at~$p$ and ending at a common endpoint, of length
$L/2$.
\item[(ii)] Any point of self-intersection of the loop~$\gamma$ is no
further than~$\tfrac{1}{2} \left(L - \sys(X) \right)$ from~$p$.
\end{enumerate}
\end{lemma}

\begin{proof}
Consider an arclength parametrisation~$\gamma(s)$ with
$\gamma(0)=\gamma(L)=p$.  Let~$q = \gamma \left( \tfrac{L}{2} \right)
\in X$ be the ``midpoint'' of~$\gamma$.  Then~$q$ splits~$\gamma$ into
a pair of paths of the same length~$\tfrac{L}{2}$, joining~$p$ to~$q$.
By Remark~\ref{ball2}, if~$q$ were contained in the open ball
$B(p,\tfrac{L}{2})$, the loop~$\gamma$ would be contractible.  This
proves item~(i).

\forget Let~$c$ be a minimizing path from~$x$ to~$x'$, and form the
loops~$\gamma_{1} \cup c$ and~$c\cup\gamma_{2}$.  Since the
loop~$\gamma$ is noncontractible, so is one of the new loops,
say~$\gamma_{1} \cup c$.  Since it is based at~$x$, we have
$\length(\gamma_{1} \cup c) \geq L$.  Therefore,
\[
L \leq \length(\gamma_{1}) + \length(c) \leq 2 \length(\gamma_{1}) L.
\]
Thus,~$c$ and~$\gamma_{1}$ have the same length, so that~$\gamma_{1}$
is length-minimizing.  The same holds for~$\gamma_{2}$, 
\forgotten

If~$p'$ is a point of self-intersection of~$\gamma$, the loop~$\gamma$
decomposes into two loops~$\gamma_{1}$ and~$\gamma_{2}$ based at~$p'$,
with~$p\in \gamma_{1}$.  Since the loop~$\gamma_{1}$ is shorter than
the pointed systolic loop~$\gamma$ at~$p$, it must be contractible.
Hence~$\gamma_{2}$ is noncontractible, so that
\[
\length(\gamma_{2}) \geq \sys(X).
\]
Therefore,
\[
\length(\gamma_{1}) = L - \length(\gamma_2) \leq L - \sys(X),
\]
proving item (ii).
\end{proof}

\begin{definition}\label{loose}
Let~$\gamma$ be a pointed systolic loop at~$p \in X$, and let
$L=\length(\gamma)=\sys(X,p)$. Let~$r$ be a real number satisfying 
\[
L - \sys(X) <2r< L.
\]
If for at least one such~$r$, the loop~$\gamma$ meets two distinct
connected components of the level curve~$S(r)=\{x\in
X|\dist(x,p)=r\}\subset X$, then we say that~$\gamma$ is a {\em loose
loop\/}.
\end{definition}

\begin{prop}\label{prop:split}
The loop homotopy class~$[\gamma]$ of a loose loop~$\gamma$ splits off in 
$\pi_1(X,p)$, \cf Definition~$\ref{def:split}$. In particular,~$\FIG(X)>0$.
\end{prop}

In other words, there is a 2-complex~$Y$ and a map~$X\to Y \vee S^1$
inducing an isomorphism~$\varphi :\pi_1(X,p) \to \pi_1(Y, y_0) * \Z$
where the image of~$[\gamma]$ is a generator of the free
factor:~$\varphi ([\gamma]) = 1\in \pi_1(S^1)=\Z$.

\begin{proof}
Consider the natural projection~$h:X \to \Reeb(f,r)$, where~$r$ is as in 
\defref{loose} so that~$\gamma$ meets two different 
components of~$S(r)$.
By Lemma~\ref{lem:simple}, the image~$h(\gamma)$ of~$\gamma$ is homotopic to the 
simple loop~$c=\alpha \cup \beta$ where~$\alpha$ agrees with the unique embedded 
arc of~$T_{r}$ with the same endpoints as~$\beta = \gamma \cap [f>r]$.
Therefore, the loop homotopy class of~$h(\gamma)$ splits off 
in~$\pi_{1}(\Reeb(f,r))$.
This yields the desired result since the natural
projection~$h$ induces an isomorphism of fundamental groups by
\lemref{24}.
\end{proof}

\section{Loose loops {\em vs} area lower bounds}
\label{three}

The following proposition provides a lower bound for the length of
level curves in a~$2$-complex{}~$X$.

\begin{proposition} 
\label{lem:1compo}
Let~$L=\sys(X,p)$ be the pointed systole of a finite~$2$-complex.
Let~$r$ be a real number satisfying
\[
L - \sys(X) <2r< L .
\]
Consider the level curve~$S \subset X$ at distance~$r$ from~$p\in X$.
Let~$\gamma$ be a pointed systolic loop at~$p$.  If~$\gamma$ is not
loose, then
\begin{equation} \label{eq:lengthS}
\length S \geq 2r - L + \sys(X).
\end{equation}
\end{proposition}

\begin{proof}
By Lemma~\ref{lem:simple}, the loop~$\gamma$ is formed of two
distance-minimizing arcs which do not meet at distance~$r$ from~$p$.
Thus, the loop~$\gamma$ intersects~$S$ at exactly two points.
Let~$\gamma' = \gamma \cap B$ be the subarc of~$\gamma$ lying
in~$B=B(p,r)$.

Since~$\gamma$ meets exactly one connected component of~$S$, there exists
an embedded arc~$\alpha \subset S$ connecting the endpoints of
$\gamma'$.  By Remark~\ref{ball2}, every loop based at~$p$ and lying
in~$B(p,r)$ is contractible in~$X$.  Hence~$\gamma'$ and~$\alpha$ are
homotopic (keeping endpoints fixed), and the loop~$\alpha \cup (\gamma
\setminus \gamma')$ is homotopic to~$\gamma$.  Hence,
\begin{equation}
\length(\alpha) + \length(\gamma) - \length(\gamma') \geq \sys(X).
\end{equation}
Meanwhile,~$\length(\gamma) =L$ and~$\length(\gamma')=2r$, proving the
lower bound \eqref{eq:lengthS}, since~$\length(S) \geq
\length(\alpha)$. 
\end{proof}

\forget
\begin{lemma} 
\label{lem:split}
In the notation of Proposition~\ref{lem:1compo}, if~$\gamma$ meets a
pair of distinct connected components of~$S$, then~$\gamma$ splits
off.
\end{lemma}

\begin{proof}
section~\ref{sec:grush}.  We use the same notation as in the proof of
Lemma~\ref{42}, where~$f$ is the distance function from~$p$, and
$C_{0}$ and~$C_{1}$ are the two connected components of~$S$ joined by
the embedded arc~$\gamma \cap [f\geq r]$ of~$Z$.  The image of this
arc by the quotient map~$Z \longrightarrow Y$ is an embedded
arc~$\alpha$ joining two points of~$A$.  By the second part of
Lemma~\ref{wedge}, the image~$\overline{\alpha}$ of~$\alpha$ under the
quotient map~$Y \longrightarrow Y/A$ splits off.  This image agrees
with the one of~$\gamma$ under the quotient map~$q:X \longrightarrow
X/B$, where~$Y/A \subset X/B$ and~$B=[f \leq r]$.  By
Remark~\ref{ball2} and Lemma~\ref{pi-exc}, the map~$q$ induces an
isomorphism of fundamental groups.  Therefore, the loop~$\gamma$
splits off, as well.
\end{proof}

\begin{proof}[Alternative proof]
Consider the Reeb graph~$T_r$ of the ball~$B=B(p,r)$ alone (rather
than all of~$X$).  Then~$T_r$ is a tree.  Consider the hybrid quotient
space~$\Reeb_r$ of~$X$, by only collapsing the connected components up
to distance~$r$.  We have~$T_r\subset \Reeb_r$.  Let~$h: X\to \Reeb_r$
be the projection.  Here the complement~$X \setminus B$ maps
homeomorphically into~$\Reeb_r$.  The map~$h$ induces isomorphism of
fundamental groups.

If~$\gamma \cap S$ has two points,~$s$ and~$t$, lying in different
connected components of~$S$, we consider the unique embedded path
$\gamma_r\subset \Reeb_r$ from~$h(s)$ to~$h(t)$ in the tree~$T_r
\subset \Reeb_r$.  The path defines a relative class~$[\gamma_r]$ for
the pair~$(\Reeb_r, X\setminus B)$.  If we choose the ``basepoint"
in~$\Reeb_r$ to be the path~$\gamma \setminus B$, then the
image~$[\gamma_r] \in \pi_1(\Reeb_r, \gamma\setminus B)$ of~$[\gamma]$
under~$h_r$ is the class that splits off.
\end{proof}
\forgotten

\section{The bound~$\SR(X) \leq 4$ for Grushko unfree complexes}
\label{sec:zero}

The length estimate of the previous section easily yields a uniform
bound for the systolic ratio of any Grushko unfree 2-complex, \ie a
complex{}~$X$ with~$\FIG(X)=0$.

\begin{theorem} 
\label{theo:smallballs}
Let~$X$ be a finite piecewise flat~$2$-complex.  Let~$\gamma
\subset X$ be a pointed systolic loop at~$x \in X$, and assume
$\gamma$ is not loose.  For every real number~$r$ such that
\begin{equation}
\label{51z}
\sys(X,x)-\sys(X) \leq 2r \leq \sys(X,x),
\end{equation}
the area of the ball~$B(x,r)$ of radius~$r$ centered at~$x$ satisfies
\begin{equation} 
\label{eq:smallballs}
\area B(x,r) \geq \left( r- \tfrac{1}{2}(\sys(X,x)-\sys(X))
\right)^{2}.
\end{equation}
\end{theorem}

\begin{proof}[Proof of Theorem~$\ref{theo:smallballs}$] 

Let~$L= \sys(X,x)$.  Denote by~$S=S(r)$ and~$B=B(r)$, respectively,
the level curve and the ball of radius~$r$ satisfying~$L-\sys(X)<2r <
L$, centered at~$x \in X$.

The non-loose loop~$\gamma$ meets a single connected component
of~$S$.  Let~$\eps= L - \sys(X)$. Now, Proposition~\ref{lem:1compo}
implies that
\begin{equation}
\length S(r) \geq 2r - L + \sys(X) =2r-\eps.
\end{equation}
Following~\cite[5.1.B]{Gr1} and~\cite{heb}, we use the coarea formula,
\cite[3.2.11]{Fe1}, \cite[p.~267]{Ch2} to obtain
\begin{eqnarray*}
\area B(x,r) & \geq & \int_{\tfrac{\eps}{2}}^{r} \length S(\rho) \,
d\rho \\ & \geq & \int_{\tfrac{\eps}{2}}^{r} (2 \rho - \varepsilon) \,
d\rho \\ & \geq & \left( r-\tfrac{\eps}{2}\right)^{2}
\end{eqnarray*}
for every~$r$ satisfying \eqref{51z}.
\end{proof}

\begin{corollary}
\label{63}
If~$X$ admits a systolic loop which is not loose, then~$\SR(X) \leq
4$.
\end{corollary}

\begin{proof}
If~$\gamma$ is a systolic loop of~$X$, then~$\varepsilon=0$.  By
hypothesis,~$\gamma$ is not loose, and we let~$r$ tend to
$\tfrac{1}{2} \sys(X)$.
\end{proof}

\begin{corollary}
\label{theo:balls}
Let~$X$ be a Grushko unfree 2-complex.  Then the area of every
ball~$B(x,r) \subset X$ of radius~$r<\frac{1}{2} \sys(X)$ centered at
a point~$x$ lying on a systolic loop of~$X$, satisfies the lower bound
\begin{equation} \label{eq:small}
\area B(x,r) \geq r^{2}.
\end{equation}
In particular, we have the bound~$\SR(X)\leq 4$.
\end{corollary}

\begin{example}
The Moore spaces~$M_n$ with~$\pi_1(M_n) = \Z_n$ satisfy the bound
$\SR(X)\leq 4$ for all~$n$.
\end{example}

\section{Noncontractible fibers} 
\label{sec:fibers}

Corollary~\ref{theo:balls} is limited by its hypothesis on the free
index of Grushko.  As a warm-up to the more exotic uses of the coarea
formula in later sections, we first present a simple argument in
Section~\ref{corank}, to obtain an explicit upper bound for~$\SR(X)$
dependent on the corank of the fundamental group.  The topological
result of this section will be used in Section~\ref{corank}.

\forget
\begin{proposition}
\label{mainloop}
Let~$Y$ be a finite connected~$2$-dimensional complex, and consider
a distance function ~$Y \longrightarrow \R$ whose Reeb graph is a
tree.  If every fiber of the lifted map~$f$ is a connected graph which
is contractible in~$Y$, then~$Y$ is simply connected.
\end{proposition}
\forgotten

\begin{lemma}
\label{halfinterval}
Let~$h: Y \to [0,1]$ be a map that yields a trivial bundle over the
half-interval~$(0,1]$. Let~$F=h^{-1}(0)$ and assume that~$F$ is
compact. Assume also that the pair~$(Y,F)$ is~$($homeomorphic to$)$ a
$CW$-pair. If every fiber of~$h$ is contractible in~$Y$ and path
connected then~$Y$ is simply connected.
\end{lemma}

\begin{proof} Since ~$(Y,F)$ is a~$CW$-pair, there exists a neighborhood~$U$ of  
$F$ such that~$F$ is a deformation retract of~$U$. Since~$F$ is compact, there 
exists~$\delta>0$ such that~$h^{-1}[0,\delta)\subset U$.  Let 
$A=h^{-1}[0,\delta)$,~$B=h^{-1}(0,1]$ and~$C=A\cap B$. The homomorphism 
$\pi_1(C) \to \pi_1(B)$ (induced by the inclusion) is an isomorphism, while the 
homomorphism~$\pi_1(C) \to \pi_1(A)$ is trivial since every loop in~$C$ is 
homotopic in~$A$ to a loop in ~$U$, and therefore to a loop in~$F$. So, by the 
Seifert -- van Kampen Theorem,~$Y$ 
is simply connected.
\end{proof}

\begin{cory}
\label{tree}
Let~$T$ be a tree and~$h: Y \to T$ be a map that yields a trivial bundle over 
the interior of each of edges of~$T$. Assume also that for every~$t\in T$ the 
pair~$(Y, h^{-1}(t))$ is a~$CW$-pair. Finally, assume that~$h^{-1}(T')$ is 
compact for every finite subtree~$T'$ of~$T$. If every fiber of~$h$ is 
contractible in~$Y$ and path connected then~$Y$ is simply connected.
\end{cory}

\begin{proof}
We first consider the case of~$T=[0,1]$. Then~$Y=A\cup B$ with
$A=h^{-1}[0,1)$ and~$B=h^{-1}(0,1]$. Notice that all the fibers of~$h$
are compact. So,~$A$ and~$B$ are simply connected by
\lemref{halfinterval}. Now the result follows from the Seifert -van
Kampen Theorem.

Now consider the case of a finite tree~$T$.  The result follows from the
Seifert -van Kampen Theorem by an obvious induction.

In the general case, consider a loop~$\gf: S^1\to Y$.  The image
of~$h\gf$, which is compact, is contained in some finite subtree~$T'$
of~$T$. Now let~$Y'=h^{-1}(T')$.  By Step 2, the loop~$\gf$ is
contractible in~$Y'$. Hence the result.
\end{proof}

\forget
\begin{proposition}If every fiber of~$h$ is contractible in~$Y$ and path 
connected then~$Y$ is simply connected.
\label{mainloop}
Let~$Y$ be a connected~$2$-dimensional complex~$($not necessarily
finite$)$, and let~$T$ be a tree. Let~$h:Y \longrightarrow T$ be a map
that yields a trivial bundle over the interior of each of edges
of~$T$.  If every fiber of~$h$ is a connected finite graph which is
contractible in~$Y$, then~$Y$ is simply connected.
\end{proposition}

\begin{proof}
First, consider the case of finite complex~$Y$.  Without loss of
generality, we can assume that the map~$f: Y\to T$ to the tree is
onto.  Now a leaf~$A \in T$ by definition has only one neighboring
vertex~$B \in T$.  Consider the edge~$[AB]\subset T$.  Let~$T'\subset
T$ be the subtree obtained by removing~$A$ and the interior of the
edge~$[AB]$.  \forgotten

\forget
Note that along the interior of the edge, we have a locally trivial bundle.  
\forgotten

\forget If the inclusion of the fiber over~$A$ in~$Y$ induces a
trivial map in~$\pi_1$, and also the inclusion of the fiber over a
typical point of the edge induces a trivial map in~$\pi_1$, then by
Seifert -- van Kampen theorem, the inclusion~$f^{-1}(T') \to Y$
induces an isomorphism in~$\pi_1$, and the proof is completed by
induction.

Now, if~$Y$ is not finite, we consider all its finite subcomplexes{} satisfying 
conditions of the Proposition, and conclude that each of them is simply 
connected. Then~$Y$ must be simply connected, since every loop in~$Y$ is 
contained in such a finite subcomplex.
\end{proof}
\forgotten

\begin{proposition} 
\label{mainloop} 
Let~$X$ be a connected~$2$-dimensional complex, and let ~$f: X \to \R$ be
a semialgebraic function.  If the group~$\pi_1(X)$ is not free, then,
for a suitable~$\rho\in \R$, the preimage~$f^{-1}(\rho)$ contains a
noncontractible loop of~$X$.
\end{proposition}

\begin{proof}  
Consider the Reeb graph~$\Reeb(f)$ of~$f$ and let~$\wh f: X \to
\Reeb(f)$ be the lifted function. Let~$p: T \to \Reeb(f)$ be the
universal cover of~$\Reeb(f)$.  Consider the pull-back diagram
\[
\CD
Y @>>> X\\
@VhVV @VV \wh fV\\
T @>p>> \Reeb(f)
\endCD
\]
Since~$\pi_1(X)$ is not free, the group~$\pi_1(Y)$ is non-trivial.
Hence, by \coryref{tree} and in view of \propref{prop:reeb}, there
exist~$t\in T$ and a loop in~$h^{-1}(t)$ that is homotopy non-trivial
in~$Y$. Therefore, there exists a loop in~$\wh f^{-1}(p(t))$ which is
non-contractible in~$X$.  This loop is contained in~$f^{-1}(\rho)$
with~$\rho=f(\wh f^{-1}(p(t)))$.
\end{proof}

\section{Corank-dependent inequalities}
\label{corank}

In this section, we apply the coarea formula to obtain an explicit
bound for~$\SR(X)$ dependent on the corank of the fundamental group.
More refined applications of the coarea formula will be presented in
later sections, to obtain a uniform bound.

We will say that a space~$X$ is {\em unfree\/} if its fundamental
group is not free, and that~$X$ is of zero corank if~$\pi_1(X)$ is.
The first theorem of this section is a special case of
Corollary~\ref{theo:balls} (but see Remark~\ref{113}).  We include its
proof as preparation for the more general result of Theorem~\ref{41z}.

\begin{theorem}
\label{41}
Every~$2$-dimensional, non-simply-connected, piecewise flat
complex~$X$ of zero corank satisfies the bound~$\SR(X)\leq 4$.
\end{theorem}

\begin{proof}[Proof of Theorem~$\ref{41}$]
Choose a point~$p\in X$ and consider the function~$f(x)=\dist(p, x)$.
The idea is to apply the coarea formula ``inward'', \ie toward the
basepoint~$p\in X$, so as to obtain the desired lower bound for the
area, along the lines of \cite[p.~129, line~9]{Gr1}.

In fact, we will use the coarea formula to show that the area of the
2-complex is bounded below by the area of a right triangle with base
$\tfrac{1}{2}\sys(X)$ and altitude~$\sys(X)$.

Consider the Reeb graph~$\Reeb(f)$ of~$f$ and note that the preimages
of~$f$ are finite graphs and~$CW$-subspaces of~$X$, \cf \theoref{46lc}
and \corref{homeo}.  By Proposition~\ref{mainloop}, there exists a
connected component
\[
C_{\rho} \subset f^{-1}(\rho) \subset X,
\]
containing a loop which is noncontractible in~$X$.

Since~$X$ is a compact path metric space, it satisfies the Hopf--Rinow
theorem,~\cite[p.~9]{Gr3}.  For each~$x\in C_\rho$, consider a
minimizing geodesic path joining~$p$ to~$x$.  Given~$r < \rho$, let
$x^r$ be a point where this geodesic meets the level set~$f^{-1}(r)$.
Then
\[
\dist(x^r, x)=f(x)-f(x^r)=\rho -r .
\]
The points~$x^{r}$ lie in the boundary of the connected component of
the superlevel set~$[f \geq r]$ containing~$C_{\rho}$.  Thus, the
connected components~$C^{r}_{1},\dots,C^{r}_{k}$ of this boundary
coalesce forward.  By Lemma~\ref{42}, there is only one such
component, denoted~$C^{r}$ (that is,~$k=1$ and~$C^{r}=C^{r}_{1}$).

\begin{lemma}
We have
\begin{eqnarray}
\rho & \geq & \tfrac{1}{2} \sys(X) \label{eq:ib} \\ \length(C^{r}) &
\geq & \sys(X) - 2(\rho-r) \label{eq:ii}
\end{eqnarray}
\end{lemma}

\forget If all the loops~$[p,x] \cup [x,y] \cup [y,p]$ were
contractible for every~$x,y \in C_{\rho}$ close enough (here,~$[p,x]$
and~$[y,p]$ are minimizing paths in~$X$ and~$[x,y]$ is a minimizing
path in~$C_{\rho}$), then the noncontractible loop contained
in~$C_\rho$ would admit a continuous retraction to~$p$, contradicting
the hypothesis.  Thus, one of these loops is not contractible.  Hence,
$2 \dist(p,C_{\rho}) = 2\rho \geq \sys(X)$.  \forgotten

\begin{proof}
The first inequality follows from Remark~\ref{ball2}.  For the second
inequality, we can assume that~$C^{r}$ does not contain a loop
noncontractible in~$X$, otherwise the result is obvious.  Then, there
is a pair of the points~$x^r, y^r\in C^r$ with~$x=y \in C_\rho \subset
X$, and a path~$[x^r, y^r]\subset C^{r}$ between them, such that the
loop~$[x^r, y^r] \cup [y^r,y]\cup [y, x^r]$ is noncontractible in~$X$.
Indeed, if all such loops were contractible, then the noncontractible
loop contained in~$C_\rho$ would admit a continuous retraction
to~$C^r$, which is impossible.  Therefore,~$\length([x^r,y^r]) \geq
\sys(X) - 2(\rho-r)$.
\end{proof}

\forget
There are two possibilities.  The set~$C^r$ may contain a loop
noncontractible in the 2-complex.  In this case, its 1-Hausdorff
measure is at least the systole of~$X$, which is enough for the
calculation below.  The remaining case is described as follows.

\begin{lemma}
Assume~$C^r$ does not contain a loop which is noncontractible in~$X$.
Then there is a pair of the points~$x^r, y^r\in C^r$ with~$x=y \in
C_\rho \subset X$, and a path~$[x^r, y^r]\subset C^{r}$ between them,
such that the loop
\[
[x^r, y^r] \cup [y^r,y]\cup [y, x^r]
\]
is noncontractible in the 2-complex, and therefore
\[
\length([x^r,y^r]) \geq \sys(X) - 2(\rho-r).
\]
\end{lemma}

\begin{proof}
If all such loops were contractible, then the noncontractible loop
contained in~$C_\rho$ would admit a continuous retraction to~$C^r$,
contradicting the hypothesis.
\end{proof}

In both cases, we have
\begin{equation}\label{diam6}
\mathcal \length(C^r)\ge \sys(X) - 2(\rho-r).
\end{equation}
\forgotten

Returning to the proof of Theorem~\ref{41}, we use the coarea formula,
\cite[3.2.11]{Fe1}, \cite[p.~267]{Ch2}, and exploit~\eqref{eq:ii} to
write
\[
\begin{aligned}
\area(X) & = \int_0^{\infty} \mathcal \length(f^{-1}(r)) \, dr \\ &
\geq \int_0^{\rho} \mathcal \length(C^r) \, dr\\ & \geq \int_{\rho
-\tfrac{1}{2}\sys(X)}^{\rho } \sys(X) -2(\rho-r) \, dr \\ & \geq
\tfrac{1}{4} \sys(X)^2,
\end{aligned}
\]
\cf \cite[p.~129, line~9]{Gr1}, proving the theorem.
\end{proof}

\begin{remark}
\label{113}
The estimate~$\SR(X)\leq 4$ can be improved in terms of the difference
$\rho-\tfrac{1}{2}\sys(X)$, where~$\rho$ is the noncontractible level.
Note that the method of proof of Corollary~\ref{theo:balls} does not
allow for such an improvement.
\end{remark}

\begin{lemma}
\label{51b}
Let~$C_\rho \subset X$ be a level curve containing a loop which is
noncontractible in~$X$.  Normalize the systole to the value~$2$.  Then
for every~$n\in \N$, the level~$C_\rho$ admits a~$\tfrac{1}
{n}$-separated set containing at least~$n+1$ elements.
\end{lemma}

\begin{proof}
Choose an essential loop~$\ell \subset C_\rho$.  Let~$A,B\in \ell$ be
a pair of points realizing the diameter of~$\ell$.
Thus~$\dist(A,B)\geq 1$.  Let~$\alpha$ be a connected component of the
complement~$\ell\setminus \{A,B\}$.  Choose a maximal finite sequence
of points~$x_i \in \alpha$,~$i=0,1,2,\ldots$ satisfying~$\dist(x_i, A)
= \tfrac{i}{n} \leq 1$.  By the triangle inequality, for all~$i \neq
j$, we have $\dist(x_i,x_j) \geq \tfrac{1}{n}$.
\end{proof}

\begin{theorem}
\label{41z}
Every unfree,~$2$-dimensional, piecewise flat complex~$X$ satisfies
the bound~$\SR(X) \leq 16 (\corank(\pi)+1)^2.$
\end{theorem}

\begin{proof}
Let~$n-1$ be the corank of~$\pi_1(X)$.  Normalize the systole to the
value 2.  By Lemma~\ref{51b}, there exists a~$\tfrac{1}{n}$-separated
set~$\{x_i\}\subset C_\rho$, such that~$\dist(x_i, x_{j}) \geq
\tfrac{1}{n}$ when~$i\not= j$.  For each~$x_i$, consider a minimizing
geodesic path joining~$p$ to~$x_i$.  Given~$r <\rho$, let~$x_i^r$ be a
point where this geodesic meets the level set~$f^{-1}(r)$.  Then for
each~$i=0,1,\ldots, n$, we have
\[
\dist(x_i^r, x_i)=f(x_i)-f(x^r_i) = \rho -r .
\]
Note that the connected set~$C_\rho$ is contained in the superlevel
set~$[ f \geq \rho]$.  Thus the set~$\{x_i\} \subset C_\rho$ coalesces
forward.  Hence by Lemma~\ref{42}, the number of components
of~$f^{-1}(r)$ which meet the set~$\{x_i\}$ is at most~$n$.  Hence
there is a pair of points~$x_k^r, x_\ell^r$ in a common connected
component of~$f^{-1}(r)$.  Let~$C_r$ be such a component.  By the
triangle inequality,
\begin{equation}
\label{61}
\begin{aligned}
\dist(x_k^r, x_\ell^r) & \geq - \dist(x_k^r, x_k) + \dist(x_k, x_\ell)
- \dist(x_\ell, x_\ell^r) \\ & \geq \tfrac{1}{n} -2(\rho-r).
\end{aligned}
\end{equation}
Therefore
\begin{equation}
\label{diam}
\length(C_r)\ge\diam C_r\ge\dist(x_k^r, x_l^r) \ge \tfrac {1}{n}
-2(\rho-r).
\end{equation}
As before, we apply the coarea formula to obtain
\[
\begin{aligned}
\area(X) & = \int_0^{\infty} \length(f^{-1}(r)) \, dr \\ & \geq
\int_0^{\infty} \length(C_r) \, dr\\ & \geq \int_{\rho -(2n)^{-1}}^{\rho
}  \tfrac{1}{n} -2(\rho-r) \, dr \\ & = \tfrac{1}{4} n^{-2},
\end{aligned}
\]
\cf \cite[p.~129, line~9]{Gr1}.  We conclude that~$4 n^2 \area(X)\ge
1$, proving the theorem.
\end{proof}

\section{Comparison with Lusternik-Schnirelmann category}
\label{two}

Originally we were led to consider the systoles of 2-complexes{} in
the context of the comparison with the Lusternik-Schnirelmann category
$\cat$, \cf \cite{KR1}. Since the latter equals 2 unless the group is
free, \cf Theorem~\ref{dim2}, the question arose whether an unfree
2-complex{} always satisfies a systolic inequality. Eventually we
found an affirmative answer in \cite{Gr1}, \cf
inequality~\eqref{eq:72}. Thus, the systolic category~$\syscat$
of~$X$, defined in \cite{KR1}, coincides with~$\cat(X)$ for every
2-complex~$X$, \cf Theorem~\ref{theo:syscat}:
\begin{equation}
\label{11}
\syscat = \cat.
\end{equation}
Note that the two categories coincide, as well, for arbitrary closed
3-manifolds \cite{KR2}. An open question is whether all Poincar\'e
3-complexes{} also satisfy equality~\eqref{11}.  Recent examples due to
J. Hillman~\cite{Hi2} show that the answer may not be easy to obtain.

\forget An~$n$-dimensional orientable manifold~$X$ is said to be {\em
essential\/} if its classifying map induces a nontrivial homomorphism
in the~$n$-dimensional integral homology group, \cf \cite{Gr1}.  If
$X$ is nonorientable, we replace~$\Z$ by~$\Z_2$ coefficients in the
definition of an essential manifold.  Clearly, an
essential~$n$-manifold is~$n$-essential.  One can extend the notion of
an essential manifold to complexes{} by mimicking the definition for
manifolds.  At any rate, the converse is not true, even for orientable
complexes~\cite[\S 6]{ba04}.  However, these two notions coincide for
closed orientable manifolds~\cite{Ba}.  \forgotten

The Lusternik-Schnirelmann category~$\cat(X)$ of a~$2$-complex~$X$
can similarly be characterized in terms of its fundamental group.
Thus, we have~$\cat X=1$ if~$\pi_1(X)$ is free, and~$\cat X=2$
otherwise (see Theorem~\ref{dim2} for a detailed statement).

\begin{theorem}
\label{dim2}
Let~$X$ be a~$2$-dimensional connected finite~$CW$-space, and let
$\pi\not =0$ denote the fundamental group of~$X$. The following
conditions are equivalent:

\begin{enumerate}
\item [(1)] the group~$\pi$ is free;

\item [(2)] the space~$X$ is homotopy equivalent to a wedge
of a
finite number of circles and~$2$-spheres;

\item [(3)] one has~$\cat X = 1$;

\item [(4)] for every group~$\tau$, every map~$f: X \to
K(\tau,1)$ can
be deformed into the~$1$-skeleton of~$K(\tau,1)$.

\end{enumerate}
\end{theorem}

\begin{remark}
Thus, a finite~$2$-dimensional complex{} satisfies a systolic
inequality if and only if none of the four equivalent conditions of
Theorem~\ref{dim2} holds.
\end{remark}

\begin{proof}
We prove the following implications:~$(1)\Rightarrow (2)\Rightarrow
(3)\Rightarrow (1)\Rightarrow (4)\Rightarrow (1)$.

The implication~$(1)\Rightarrow (2)$ is proved by C. Wall
\cite[Proposition~3.3]{Wa}.  We recall his argument for completeness.
Let~$\Z[\pi]$ denote the group ring of the group~$\pi$ that we assumed
to be free. First, Wall proved that~$\pi_2(X)$ is a finitely
generated projective~$\Z[\pi]$-module. Now, by a theorem of H.~Bass
\cite{Bas}, every finitely generated projective module over~$\Z[\pi]$
is free. Consider a wedge~$Y$ of~$k$ circles and~$l$ spheres,
where~$k$ is the number of free generators of~$\pi$ and~$l$ is a
number of free generators of the~$\Z[\pi]$-module~$\pi_2(X)$. Now we
map~$Y$ to~$X$ by mapping circles to free generators of~$\pi_1(X)$ and
spheres to~$\Z[\pi]$-free generators of~$\pi_2(X)$, and this map is a
homotopy equivalence.

The implication~$(2)\Rightarrow (3)$ is obvious. The implication
$(3)\Rightarrow (1)$ is well known, \cf \cite[Exercise 1.21]
{CLOT}. The implication~$(1)\Rightarrow (4)$ holds for~$\tau=\pi$
since~$K(\pi,1)$ is a wedge of circles if~$\pi$ is free. For
general~$\tau$, we notice that every map~$f: X \to K(\tau, 1)$ factors
through~$K(\pi,1)$ . To prove that~$(4)\Rightarrow (1)$, it suffices
to prove the implication for a map~$f$ that induces an isomorphism of
fundamental groups. Let~$K$ denote the 1-skeleton of~$K(\pi,1)$.
Then the map~$X \to K \subset K(\pi,1)$ induces an isomorphism of
fundamental groups. Thus,~$\pi$ is a subgroup of a free group
$\pi_1(K)$. Hence~$\pi$ is free.
\end{proof}

\begin{proposition}\label{homdim}
Let~$X$ be a connected, finite,~$n$-dimensional~$CW$-complex,
where
$n>2$. Let~$\pi_1(X)=\pi$ and assume that~$H^k(\pi;G)=0$ for
all
$k\ge n$ and all~$\Z[\pi]$-modules~$G$. Then~$\cat X <n$.
\end{proposition}

\begin{proof} By theorem of I. Berstein~\cite{Ber}, \cf
\cite[Theorem~2.51]{CLOT}, we have~$\cat X=n$ only if
\begin{equation}
\label{51}
u^n\ne 0 \mbox{ for some } u\in H^1(X;I(\pi)^{\otimes n}),
\end{equation}
where~$I(\pi)$ is the augmentation ideal of~$\Z[\pi]$. We can obtain
the classifying space~$K(\pi,1)$ by attaching~$k$-cells with~$k\ge 3$
to~$X$. Now, the inclusion~$X \to K(\pi,1)$ induces an isomorphism
$H^1(\pi;G) \to H^1(X;G)$ for any~$\Z[\pi]$-module~$G$.  Therefore
$u^n=0$.
\end{proof}

\begin{cory}
Every free, connected, finite~$CW$-complex~$X$ of dimension at least
$2$ satisfies~$\cat(X)\leq n-1$.
\end{cory}

\begin{proof}
Notice that~$H^i(F;G)=0$ for every free group~$F$ and~$i>1$.  Now, for
$n>2$ the claim follows from \propref{homdim}, while for~$\dim X=2$
the claim follows from Theorem~\ref{dim2}, item~(3).
\end{proof}

An invariant called {\em systolic category\/},~$\syscat(X)$, of~$X$
was defined in~\cite{KR1}. It is a homotopy invariant, which,
furthermore, coincides with the Lusternik-Schnirelmann category
$\cat(X)$ for all 3-manifolds \cite{KR2}.  We now calculate it for
2-dimensional complexes. For technical reasons, we need to describe
the 1-dimensional case first.

\begin{proposition}
Every graph~$X$ satisfies~$\syscat X = \cat X$.  The common value
is~$0$ if~$X$ is contractible, and~$1$ otherwise.
\end{proposition}

\begin{proof}
A graph~$\Gamma$ which contains nontrivial cycles, satisfies the
obvious systolic inequality~$\sys(\Gamma) \leq \length (\Gamma)$.  The
1-systole of a tree is infinite, being an infimum over an empty set. 
\end{proof}

\begin{theorem} 
\label{theo:syscat}
Let~$X$ be a~$2$-dimensional complex that is not homotopy equivalent to a
wedge of circles. Then we have~$\syscat(X) = \cat(X)$.
\end{theorem}

\begin{proof}
Every free~$X$ is homotopy equivalent to a wedge of circles and
2-spheres, \cf \cite{Wa} and Theorem~\ref{dim2}.  We replace~$X$ by
such a wedge~$W$.  For any~$K>0$, we can find a metric
with~$\sys(W)^2\geq K\area(W)$.  Clearly, the 2-systole of~$W$
satisfies~${\rm sys}_2(W)\leq \area(W)$.  Therefore~$\syscat(W)
=\cat(W)=1$ in this case. 

The theorem now follows from~\eqref{eq:72} (or Theorem~\ref{dim2}),
combined with the homotopy invariance of both categories \cite{KR2}.
\end{proof}

\section{A useful auxiliary space}
\label{sec:Z}

The main goal of this section is the construction of a space~$Z$
obtained by cutting a loose loop of~$X$ along a graph, and folding the
graph to a tree.

We assume that~$X$ is connected.  Let~$p\in X$ be a point on a
systolic loop of~$X$.  Let~$f$ be the distance function from~$p$.
Let
\[
T'_r \subset \Reeb(f,r)
\]
be the pruned Reeb tree of the ball of radius~$r< \tfrac{1}{2}
\sys(X)$, as in Section~\ref{sec:tree}.  Let~$e \subset T'_{r}$ be an
open non-root edge of length~$t_{0}$, isometrically identified
with~$(0,t_{0})$.

Let~$t \in (0,t_{0}) = e$.  Note that the connected component~$C_t
\subset X$ of the level set of~$f$ corresponding to~$t$ is
contractible in~$X$ by Remark~\ref{ball2}.  Let~$\lambda=f(C_t)=\bar
f(t) \in \R$ be the value taken by~$f$ on~$C_t$.

Next, we define new spaces~$W$ and~$Z$ as follows.  Consider the
complement~$X \setminus C_t$.  Glue back two copies of~$C=C_t$ to~$X
\setminus C$ in order to compactify the two open sets
\[
U_{-} = \pi^{-1}(e) \cap [f<\lambda] \quad \mbox{and} \quad U_{+} =
\pi^{-1}(e) \cap [f>\lambda]
\]
in the neighborhood of~$C$, where~$\pi:X \longrightarrow \Reeb(f,r)$ is
the quotient map.  Denote by~$W$ the resulting complex.  If~$C$ is
a tree, no further modifications need to be made.

\begin{remark}
\label{pl}
The space~$W$ is not precisely of the type envisioned in earlier
sections, since~$C$ is a union of circular arcs rather than straight
line intervals.  However, replacing~$C$ by a sufficiently fine
polygonal curve, we obtain a piecewise flat complex whose systole and
area differ from those of~$W$ by an arbitrarily small amount.  Thus,
for the purposes of proving our systolic bound, we may certainly
assume~$W$ is piecewise flat as defined earlier.
\end{remark}

Otherwise we introduce further identifications on~$W$ as follows.  If
the graph~$C$ contains a nontrivial embedded loop, we map this loop
isometrically to a concentric circle in the complex plane.  We then
use complex conjugation, folding the circle to an interval, to
introduce the same identification on both copies of~$C$ in~$W$.

Inductively, we can eliminate all cycles of both copies of~$C$ in~$W$
while preserving the structure of a piecewise flat complex{} (up to
subdivision).  Denote by~$C_\pm$ the pair of trees thus obtained after
elimination of all cycles.  Here~$C_\pm$ is attached to~$U_\pm$ as
before.  By construction, the new space
\[
Z= (X \setminus C) \cup C_{-} \cup C_{+}
\]
has the same fundamental group as~$\Reeb(f,r)\setminus e$, \cf
Proposition~\ref{33}.  By Remark~\ref{pl}, the space~$Z$ may be
assumed to be equipped with a piecewise flat structure induced
from~$X$, and~$\area(Z)=\area(X)$.  Let~$M(X,r)$ be the minimal model,
\cf Section~\ref{sec:tree}.

\begin{lemma}
\label{144}
Let~$e\subset M(X,r)$ be an open non-root edge.  Then the natural map
\begin{equation} 
\label{eq:i}
i:Z \longrightarrow Y= M(X,r) \setminus e
\end{equation}
sending~$U_{-} \cup C_{-}$ and~$U_{+} \cup C_{+}$ to the endpoints
of~$e$ induces an isomorphism of fundamental groups (on each connected
component).
\end{lemma}

\begin{proof}
The proof is immediate from the construction.
\end{proof}

\begin{lemma} \label{lem:FIG(Z)}
There exists an unfree connected component~$Z_{*}$ of~$Z$ with~$\FIG
(Z_{*}) \leq \FIG(X) -1$.
\end{lemma}

\begin{proof}
The lemma follows from Lemma~\ref{144} and Proposition~\ref{32}.
\end{proof}

\section{The uniform bound}
\label{fifteen}

\begin{theorem} 
\label{theo:unif}
Every finite unfree piecewise flat~$2$-complex~$X$ satisfies the bound
\[
\SR(X) \leq 12.
\]
\end{theorem}

\begin{proof}
By \corref{theo:balls}, the bound holds if~$\FIG(X)=0$.  Assume that
the bound holds if~$\FIG(X)<n$.  We will prove the bound for the case
$\FIG(X)=n$.  We use the notation of Section~\ref{sec:Z}.  Let~$e
\subset T'_r$ be an open non-root edge of the pruned tree~$T'_r
\subset M(X,r)$ where~$r < \frac{1}{2} \sys(X)$.  Consider the space
$Z_*=Z_*(e)$.  By Lemma~\ref{lem:FIG(Z)},~$\FIG(Z_*)<n$, while
$\area(Z_*)=\area(X)$.  Suppose there is a connected component~$C=C_t
\subset X$ of a level curve of the distance function, with~$t \in e$,
such that~$\sys(X) \leq \sys(Z_*)$.  Then,
\[
\SR(X)= \frac{\sys(X)^2} {\area(X)} \leq \frac {\sys(Z_*)^2}
   {\area(Z)} =\SR(Z_*)\le 12,
\]
by the inductive hypothesis.  Now let us assume that~$\sys(X) \geq
\sys(Z_*(e))$, for every non-root edge~$e \subset T'_{r}$ and all~$t
\in e$.  Then every systolic loop of~$Z_{*}$ must meet either~$C_{-}$
or~$C_{+}$.  Indeed, if a systolic loop~$\gamma\subset Z_{*}$ lies in
\[
Z \setminus \left( C_- \cup C_+ \right) = X \setminus C \subset X,
\]
then
\[
\sys(X) \leq \length(\gamma) = \sys(Z_{*}) .
\]
Denote by~$X_{\tau}$ the connected component of the level set of~$f$
corresponding to the point~$\tau \in e=(0,t_{0})$.  We have the
following lemma.

\begin{lemma} 
\label{lem:X_r} 
Let~$\gamma$ be a systolic loop of~$Z_{*}$.
\begin{enumerate}
\item If~$\gamma$ meets~$C_{-}$, then, for every~$\tau \in (0,t)$,
\[
\length(X_{\tau}) \geq 2t - 2\tau;
\]
\item If~$\gamma$ meets~$C_{+}$, then, for every~$\tau \in (t,t_{0})$,
\[
\length(X_{\tau}) \geq 2\tau - 2t.
\]
\end{enumerate}
\end{lemma}

\begin{proof}
Suppose that~$\gamma$ meets~$C_{-}$.  By Lemma~\ref{144}, the image
$i(\gamma)$ must leave the corresponding endpoint of the edge
$e\subset M(X,r)$.  In particular,~$\gamma$ must meet the level
$X_\tau$. Moreover, since~$i(\gamma)$ covers~$(0,t)$ at least twice,~$\gamma$ 
meets~$X_{\tau}$ at least twice.

Choose a subarc~$\alpha\subset U_{-} \cup C_{-}$ of~$\gamma$ which
meets~$C_{-}$, and with its endpoints in~$X_{\tau}$.  Note
that~$\length(\alpha) \geq 2(t-\tau)$.  Let~$\beta \subset X_{\tau}$
be an embedded path joining the endpoints of~$\alpha$.  By
Lemma~\ref{144}, the loop~$(\gamma \setminus \alpha) \cup \beta$ is
homotopic to the systolic loop~$\gamma$ in~$Z_{*}$.  Since it cannot
be shorter than~$\gamma$, we have
\[
\length(\beta) \geq \length(\alpha).
\]
Hence,
\[
\length(X_{\tau}) \geq 2(t-\tau),
\]
which proves item~(1).  Item~(2) is proved using similar arguments.
\end{proof}

\begin{lemma} \label{lem:last}
Assume no connected component~$($of a level curve of the distance
function$)$~$C \subset X$ can be found such that~$\sys(X) \leq
\sys(Z_*)$.  Then we have
\[
\area(\pi^{-1}(e)) \geq \tfrac{1}{2} \length(e)^{2}.
\]
\end{lemma}

\begin{proof}
Let~$A_-$ and~$A_+$ denote the set of values~$t$ for which there exists a 
systolic 
loop of~$Z_*$ that meets~$C_-$ and~$C_+$, respectively. Then~$A_-$ and~$A_+$ are 
relatively closed subsets of
$(0,t_0)$.  Hence, if~$A_+\cap A_-=\emptyset$ then at least one of the sets 
$A_-$, 
$A_+$ is
the full interval~$(0,t_{0})$.

Suppose~$A_+\cap A_-$ is nonempty.  Then for some~$t \in [0,t_0]$,
there exist two systolic loops~$\gamma_{-}$ and~$\gamma_{+}$ of~$Z_*$
such that~$\gamma_{\pm}$ meets~$C_{\pm}$.  Now items (1) and~(2) of
Lemma~\ref{lem:X_r} provide a lower bound for the length of~$X_{\tau}$
for every~$\tau \in (0,t_{0})$.  Integrating this lower bound from~$0$
to~$t_{0}$ leads, through the coarea formula, to the lower bound
\[
\area(\pi^{-1}(e)) \geq \tfrac{1}{2} t_{0}^{2}.
\]

Suppose~$A_- =(0,t_0)$.  Then for all~$t \in e$, there exists a
systolic loop of~$Z_*$ that meets~$C_{-}$.  Integrating the lower
bound provided by item~(1) of Lemma~\ref{lem:X_r} over~$[0,t_0]$, we
obtain, through the coarea formula, the bound
\[
\area(\pi^{-1}(e)) \geq t_{0}^{2}.
\]
The same bound holds if~$A_+ =(0,t_0)$.
\end{proof}

The pruned tree~$T'_{r}$ decomposes into the root edge~$e_{p}$ of
length~$\ell \geq 0$ (possibly zero) based at~$\pi(p)$, and two trees
$\Gamma$ and~$\Gamma'$, of height~$r-\ell$, attached to~$e_{p}$ at the
other endpoint:
\[
T'_r = e_{p} \cup \Gamma \cup \Gamma' .
\]
The bound of Lemma~\ref{lem:last} can be improved for the root
edge~$e_{p}$ to
\[
\area(\pi^{-1}(e_{p})) \geq \ell^{2}.
\]
Indeed, for every~$\tau \in [0,\ell]$, the level curve~$S_{\tau} =
f^{-1}(\tau)$ is connected.  Furthermore, every systolic loop
through~$p$ meets~$S_{\tau}$, \cf Lemma~\ref{lem:simple}.  The lower
bound of Proposition~\ref{lem:1compo} on the length of~$S_{\tau}$
leads, via the coarea formula, to the desired bound.

Thus, each edge of~$\Gamma$ and~$\Gamma'$ makes a contribution of one
half of its length squared to the total area of~$X$, while the root
edge~$e_{p}$ makes a contribution equal to the square of its length.  Hence,
\begin{eqnarray*}
\area(X) & \geq & \ell^{2} + \tfrac{1}{2} E(\Gamma) + \tfrac{1}{2}
         E(\Gamma') \\ & \geq & \ell^{2} + \tfrac{1}{2} (r-\ell)^{2}
         \\ & \geq & \tfrac{1}{3} r^{2},
\end{eqnarray*}
where the second inequality comes from Proposition~\ref{21}.  The proof
of Theorem~\ref{theo:unif} is then completed by letting~$r$ tend
to~$\tfrac{1}{2} \sys(X)$.
\end{proof}

\section{Acknowledgment}

We are grateful to Philip Boyland and Alex Dranishnikov for useful
discussions.


\end{document}